\documentclass[12pt]{article}
\usepackage{graphicx}
\usepackage{a4wide}
\usepackage[english]{babel}
\usepackage{amsfonts}
\usepackage{amsmath}
\usepackage{amssymb}
\usepackage{dsfont}

\allowdisplaybreaks

\newcommand {\ve} {\varepsilon}

\newenvironment{mypar}[2]
  {\begin{list}{}%
    {\setlength\leftmargin{#1}
    \setlength\rightmargin{#2}}
    \item[]}
  {\end{list}}

\begin{document}
\title{\bf Paul L\'evy, strong approximation and the St.\ Petersburg paradox}
\author{Istv\'an Berkes
\footnote{ Graz University of Technology, Institute of Statistics,
Kopernikusgasse 24, 8010 Graz, Austria.  \mbox{e-mail}:
\texttt{berkes@tugraz.at}. Research supported by FWF grant
P24302-N18 and NKFIH grant K 108615.} }
\date{}
\maketitle

\abstract{The St.\ Petersburg paradox, formulated in the early 1700's, concerns the 'fair'
entry fee in a game where the winnings are distributed as $P(X=2^k)=2^{-k}$,
$k=1, 2, \ldots$.
Since the tails of $X$ are not regularly varying, the accumulated gain $S_n$ in $n$ St.\ Petersburg
games has no limit distribution after any centering and norming, making the asymptotic study of
the game a challenging problem. The problem was solved by Martin-L\"of
(1985) and Cs\"org\H{o} and Dodunekova (1991), leading to a clarification of the paradox
and a fascinating  asymptotic theory. The purpose of this paper is to discuss
the crucial, but forgotten contribution of Paul L\'evy (1935) to the field. In a remark
in his classical paper \cite{levy35}, L\'evy determines the asymptotic distribution of a large class
of i.i.d.\ sums with generalized St.\ Petersburg tails $c x^{-\alpha}\psi(\log x)$,
where $0<\alpha<2$ and $\psi$ is a periodic function on $\mathbb R$.
His proof uses a coupling argument
similar to Skorohod representation and provides a strong (pointwise) approximation result, the first
in probability theory.
The argument also yields a strong approximation approach to
to domains of attraction and partial attraction, as well as strong
approximation of i.i.d.\ sums with infinitely divisible variables.

We finally discuss an argument of L\'evy \cite{levy35}, also of considerable historical interest, proving a limit
theorem  via the quantile transform, another 'first' in probability theory.
His argument yields a qualitative version of the stable decomposition theorem of  LePage,
Woodroofe and Zinn (1981) and, adapted to semistable variables, leads to a complete solution of the strong approximation
problem for St.\ Petersburg sums, as we will show in a subsequent paper \cite{buj}.

\section{Introduction}

Let $X, X_1, X_2, \ldots$ be i.i.d.\ r.v.'s
with
\begin{equation}\label{def}
P(X = 2^k) = 2^{-k}, \quad  (k = 1,2, \dots)
\end{equation}
and let $S_n=\sum_{k=1}^n X_k$. The asymptotic behavior of the
sequence $\{S_n, \, n\ge 1\}$ has  attracted considerable
attention in the literature in connection with the St.\ Petersburg
paradox (for the 'standard' formulation, see Daniel Bernoulli
\cite{bern}),  concerning the 'fair' entry fee in a game where the
winnings are distributed as $X$. We refer to Cs\"org\H{o} and
Simons \cite{cssi} for a historical account and bibliography of
the problem. Solving the entry fee problem requires determining
the precise asymptotic behavior of $S_n$.
Feller \cite{fe} proved that
\begin{equation} \label{feller1945}
\lim_{n\to\infty} \frac{S_n}{n\log_2 n}=1 \qquad \text{in
probability}
\end{equation}
(where $\log_2$ denotes logarithm with base 2) and Martin-L\"of
\cite{ml} obtained
$$S_{2^k}/2^k -k\overset{d}{\longrightarrow} G,$$
where  $G$ is the infinitely
divisible distribution function with
characteristic function $\exp (g(t))$, where
\begin{equation}\label{gdef}
g(t)=\sum^0_{l = -\infty} (e^{it 2^l} - 1 - it 2^l) 2^{-l} +
\sum^\infty_{l = 1} (e^{it 2^l} - 1) 2^{-l}.
\end{equation}
He also proved that if $n_k \sim \gamma 2^k$, $1\le \gamma <2$,
then
\begin{equation} \label{csd0}
S_{n_k}/n_k-\log_2 n_k \overset{d}{\longrightarrow} G_\gamma
\end{equation}
where $G_\gamma$ denotes the distribution with characteristic  function
$\exp(\gamma g(t/\gamma) -it \log_2 \gamma)$. Letting $\gamma_n= n/
2^{[\log_2 n]}$ (where $[ \ ]$ denotes integral part),
Cs\"org\H{o} and Dodunekova \cite{csd} proved that (\ref{csd0}) holds iff $\gamma_{n_k}\to\gamma$
and Cs\"org\H{o} (\cite{csrates}, Theorem 1)  proved that
\begin{equation}\label{csodo}
\sup_x \left | P\left( \frac{S_n}{n} - \log_2 n\le x\right) -G_{\gamma_n} (x)\right|
\longrightarrow 0 \quad \text{as} \ n\to\infty
\end{equation}
and determined the precise convergence rate. Relation
(\ref{csodo}) shows that the class of subsequential limit
distributions of $S_n/n-\log_2 n$ is the class
$$\mathcal{G}=\{ G_\gamma: 1\le\gamma<2\}.$$
If $n$ runs through the interval $[2^k, 2^{k+1}]$, then
$G_{\gamma_n}$ moves through the  distributions $G_{j/2^k}, 2^k\le
j\le 2^{k+1}$ representing, in view of $G_1=G_2$ (cf.\ \cite{ml},
Theorem 2), a "circular" path in $\mathcal G$. In view of
(\ref{csodo}), the distribution of $S_n/n-\log_2 n$ also describes
approximately a circular path, a remarkable asymptotic behavior
called {\it merging} in \cite{csrates}. Using strong approximation
of the uniform empirical process, Cs\"org\H{o} and Dodunekova
\cite{csd} showed that merging holds for extremal and trimmed sums
of the sequence $(X_n)$ as well and Berkes, Horv\'ath and Schauer
\cite{bhs} and del Barrio, Janssen and Pauly \cite{ba-ja-pa}
proved that the same holds for bootstrapped sums.

The purpose of this paper is to discuss the crucial, but forgotten
contribution of Paul L\'evy to  the field. In Section 2 we will
discuss a remark in L\'evy \cite{levy35} stating and proving the basic
limit theorem for St.\ Petersburg sums for a much larger class of
i.i.d.\ sequences.
L\'evy's proof, depending on his construction of semistable laws
and a remarkable coupling argument, provides a particularly simple
and elementary approach to the problem and fills a crucial missing
piece in St.\ Petersburg history. His method works for general
i.i.d.\ sums as well, providing a strong approximation approach to
domains of attraction and domains of partial attraction, as well
as approximation of i.i.d.\ sums with infinitely divisible
variables. Finally, we will discuss an argument in \cite{levy35}
concerning stable distributions which appears to be the first
application of the quantile transform method in probability theory
and which, adapted to semistable variables, leads to a
a complete solution of the strong approximation
problem for St.\ Petersburg sums, as we will show in a subsequent paper \cite{buj}.

\section{Asymptotics by Poisson coupling}

Let $X, X_1, X_2, \ldots$ be i.i.d.\ random variables and let $S_n=\sum_{k=1}^n X_k$. If for some numerical sequences
$(a_n)$, $(b_n)$ we have
\begin{equation}\label{dattr}
(S_n-a_n)/b_n \overset{d}{\longrightarrow} Z
\end{equation}
with a nondegenerate $Z$, then $Z$ is either Gaussian or $\alpha$-stable with some $0<\alpha<2$.
A necessary and sufficient criterion for a Gaussian limit is
\begin{equation}\label{gauss}
\lim_{x\to\infty} \frac{ x^2 P(|X|> x)}{EX^2 I\{|X|\le x\}}=0
\end{equation}
and the corresponding criterion for an $\alpha$-stable limit is that $P(|X|>x)$ is regularly varying
with exponent $-\alpha$ and
\begin{equation}
\lim_{x\to\infty} \frac{P(X>x)}{P(|X|>x)}= p, \quad  \lim_{x\to\infty} \frac{P(X\le -x)}{P(|X|> x)}= q
\end{equation}
for some $p, q \ge 0$, $p+q=1$. The criterion in the Gaussian case
is due to L\'evy (\cite{levy35}, Th\'eor\`eme II, p.\ 366 and \cite{levy37}, Th\'eor\`eme 36,3, p.\ 113); the
stable case was obtained independently by Gnedenko \cite{gn0} and
Doeblin \cite{doe1940}.  Note that their criterion was anticipated
by L\'evy, who in \cite{levy35}, p.\ 374,  remarks that

\medskip
\begin{mypar}{0.7cm}{0.7cm}
{\it I would like to stress here the profound difference between convergence to the Gaussian law [...]
and convergence to other stable laws, which is the consequence of precise hypotheses on the probability
of large values of the variable. One gets one of these stable laws as a limit only if the original law resembles that
sufficiently. [...]}
\end{mypar}

\bigskip\noindent
Then he goes on to say that

\medskip
\begin{mypar}{0.7cm}{0.7cm}
{\it There is an analogous difference between the convergence to a
Gaussian law in the case when ${\cal E}\{x^2\}$ is infinite and the case of
convergence to another stable law. Assume that all the $x_n$ have
(except perhaps in a finite interval) the same distribution which
we assume, to fix the ideas, to be symmetric and such that
$F(X)=X^{-\alpha}$ $(0<\alpha\le 2)$.  Divide the interval $(1,
\infty)$ into infinitely many intervals separated by the numbers
$X_n=q^n$ $(q>1)$, and in each of in these intervals perform an
arbitrary change of the probability distribution. If $\alpha=2$,
these modifications have no effect on the type of the limit
distribution, which is Gaussian. [...] If, on the contrary,
$\alpha<2$, then the modifications in the different intervals
cause perturbations which act successively when  $n$ grows and for
each of them there comes a moment when its effect is not any more
negligible.[...] If the modifications occur periodically, we will have
convergence to a class of associated semistable laws (here we mean
a semistable law and  its powers whose characteristic functions
are the powers of the original characteristic function).}
\end{mypar}

\medskip\noindent


A central issue
in \cite{levy35} is the connection between the behavior of
a distribution function at $\pm \infty$ and the asymptotic behavior of the
partial sums of the corresponding i.i.d.\ sequence. In addition to the main,
and now classical, theorems of the paper, L\'evy also discusses many interesting
irregular situations, one of which is the above example.
He gives no explicit direct proof of his claim, but
a proof can be recovered from  his construction of semistable
laws and a remarkable coupling argument in the case of stable
attraction.
The symmetry of the $x_k$ is not used in the proof; what L\'evy's
argument actually yields  is  that
if $0<\alpha<2$, $q \ge 2$ is an integer and $x_n$ are
i.i.d.\ random variables with tails of the form
\begin{equation}\label{tail2}
P(x_n>x) =c_1x^{-\alpha} \psi (\log_q x),\quad  P(x_n\le -x) =c_2
x^{-\alpha} \psi (\log_q x) \qquad  (x \ge x_0)
\end{equation}
where $\psi$ is a bounded function with period $1/\alpha$, then
the  class of weak subsequential limits of $n^{-1/\alpha}
\sum_{k=1}^n x_k$, suitably centered, is the class of convolution
powers of the infinitely divisible semistable distribution
whose Poisson measure has the tails in (\ref{tail2}) for
$x>0$. The periodic tail condition  (\ref{tail2}), formulated
in L\'evy's remark only verbally, is a consequence of the
structure of semistable laws, see p.\ 357, $3^{o}$.
Except a minor technical point (see Footnote
1), L\'evy's proof is complete in the case $\alpha\ne 1$; in the
case $\alpha=1$ his argument is sketchy, in particular, he does
not determine the centering factor in the limit theorem. This gap
can be removed easily by truncation,
leading to the centering factor $c\log n$.

L\'evy's argument, written in his characteristic style, is a most
interesting reading both for probabilists and anyone interested in
the history of probability theory. Below we give his proof, using
today's terminology, in the case of positive variables $x_n$.
Assume first $0< \alpha<1$, when no centering factor is needed and
the argument is the simplest. In \cite{levy35}, pp.\ 352--356
L\'evy gives his now classical construction of stable laws by
inhomogeneous Poisson processes and on pp.\ 358--359 he uses it to
prove that if $x_n$ are i.i.d.\ positive random variables with
tails $P(x_n>x)=x^{-\alpha}$ for $x\ge 1$, then $n^{-1/\alpha}
\sum_{k=1}^n x_k$ converges weakly to a completely asymmetric
stable law $\mathcal{L}$ with parameter $\alpha$. (Of course, this
holds under more general conditions and is easily proved by
characteristic functions, but L\'evy's direct coupling argument
has great historical interest.) Let $y_1>y_2> \ldots$ be the
points of an inhomogeneous Poisson process $\mathcal{P}$ on $(0,
\infty)$ with intensity measure having tails
\begin{equation}\label{tail0}
T(x)=x^{-\alpha} \qquad (x>0).
\end{equation}
Let $S=\sum y_i$, and let $\mathcal{L}$ denote the distribution of
$S$. The relation $S<\infty$ a.s.\ follows from the fact that the
expected number of points $y_i$ in $(1, \infty)$ is $T(1) <\infty$
and $E\sum y_i I(y_i\le 1) = \int_0^1 xT(dx)<\infty$ by
$\alpha<1$. Since  $nT(n^{1/\alpha} x)=T(x)$ for any $n\ge 1$, the
sum of $n$ random variables with distribution $\mathcal L$,
divided by $n^{1/\alpha}$ has again distribution $\mathcal L$ and
thus $\mathcal L$ is stable. (In accordance with his earlier book \cite{levy25},
L\'evy uses the notation $\bf{L}_{\alpha, -1}$ for this
distribution, which is the $\alpha$-stable law
with skewness parameter $\beta=-1$, see also
Gnedenko and Kolmogorov \cite{gk}, p.\ 164).
Clearly,
$${\cal L}= \lim_{t\to 0} {\cal L}^{(t)}$$
in law, where ${\cal L}^{(t)}$, $t>0$ is the distribution of the sum
$\bar{S}_t$ of the points of the Poisson process $\cal P$
exceeding $t$. The number of such points is Poisson distributed
with mean $T(t)$ and the distribution of $\bar{S}_t$ equals the
distribution of the sum of $N$ i.i.d.\ random variables
concentrated on $(t, \infty)$ with tails proportional to $T(x)$,
where $N$ is a Poisson variable having mean $T(t)$ and independent
of the i.i.d.\ sequence. Let now $x_1, x_2, \ldots$ be i.i.d.\
random variables concentrated on $(1, \infty)$ with tails $T(x)$
and $\hat{S}=n^{-1/\alpha} \sum_{k=1}^n x_k$. Clearly, $\hat{S}$
is obtained from $\bar{S}_t$ by choosing $t=n^{-1/\alpha}$ and
replacing the Poisson variable $N$ by $n$. Since the mean and
variance of $N$ is $T(n^{-1/\alpha})=n$, Chebysev's inequality
yields
$|N-n|= O_P(\sqrt{n})$, which implies
easily\footnote{This requires a probability estimate for the
maximal fluctuation of the partial sum process $\{S_j, \, j\ge 1\}$
of the sequence $(x_k)$ for $n-Cn^{1/2}\le j\le n+Cn^{1/2}$, which follows from maximal
inequalities known at the time.} that the L\'evy distance between
the distributions of $\bar{S}_t$ for $t=n^{-1/\alpha}$ and $\hat{S}$ tends to 0 as
$n\to\infty$ and, consequently, $\hat{S}
\overset{d}{\longrightarrow} {\cal L}$ as $n\to\infty$, i.e.
$n^{-1/\alpha} \sum_{k=1}^n x_k \overset{d}{\longrightarrow}  {\cal L}$.

In the case $1<\alpha<2$ the argument is similar, only from the
$x_\nu$ one has to subtract their means $Ex_\nu$ and the
definition of $S$ has to be replaced by $S=\lim_{t\to 0}
(\bar{S}_t -E\bar{S}_t)$. Here $E\bar{S}_t= \int_t^\infty
xT(dx)<\infty$ by $\alpha>1$ and the existence and a.s.\
finiteness of $S$ follows, as L\'evy remarks, from the Kolmogorov
two series criterion. More changes are needed in the case
$\alpha=1$
and L\'evy only sketches them; in particular, he does not compute the
centering factor in the limit theorem for $\sum_{k=1}^n x_k$. In
the stable case $\psi=\text{const}$ he computes the value $c_n=c
\log n$ in \cite{levy37}, p.\ 209; the argument
there works in the semistable case as well.

\medskip
Let now $\psi(x)$ be a bounded periodic function on $\mathbb R$
with period $1/\alpha$  and assume that
\begin{equation} \label{tail1}
T(x)=x^{-\alpha} \psi (\log_q x), \quad x>0
\end{equation}
is nonincreasing. Then, as L\'evy points out (see \cite{levy35},
p.\ 357, {$\rm{3^{o}}$}, "On ne retrouve...") for the special values $n=q^k$ we still
have $nT(n^{1/\alpha} x)=T(x)$ and thus replacing the Poisson
process corresponding to (\ref{tail0}) by the one corresponding to
(\ref{tail1}) and denoting the distribution of $S$ again by
$\mathcal{L}$, the sum of $n=q^k$ random variables with
distribution $\mathcal L$, divided by $n^{1/\alpha}$ has
distribution $\mathcal L$ and thus $\mathcal L$ is semistable.
Also, the previous argument yields $\hat{S}
\overset{d}{\longrightarrow} {\cal L}$ as $n\to\infty$ along these
special $n$'s, proving L\'evy's claim for the indices $n=q^k$. To
prove  the general case, in \cite{levy35}, p.\ 380 L\'evy points
out that

\medskip
\begin{mypar} {0.7cm}{0.7cm}
{\it Indeed, if for a sequence $n_p$ of $n$'s one gets laws whose
types converge to that of a law $\cal{L}'$, then for values $n_p'$
such that $n_p'/n_p$ has a limit $k$, one gets laws whose types
converge to that of ${\cal {L}}^{'k}$, which denotes the law whose
characteristic function is obtained by raising that of
${\cal{L}}'$ to the power $k$.}
\end{mypar}

\medskip\noindent
Note that this principle holds in full generality: if along a
subsequence $(n_p)$ the centered and normed partial sums of an
i.i.d.\ sequence converge  weakly to a distribution $G$, then
along another subsequence $(n_p')$ with $n_p'/n_p\to k$ where $k$
is an arbitrary positive number, the suitably centered and normed
partial sums converge weakly to the convolution power $G^{*k}$
(that is, to the law with characteristic function $\varphi^k$,
where $\varphi$ is the characteristic function of $G$). The proof
(which L\'evy omits) is immediate by characteristic functions.
Since $\varphi^k$ is the pointwise limit of a sequence of
characteristic functions, it is itself a characteristic function.
Thus in the case of the sequence $(x_n)$ in (\ref{tail2}), the
centered and normed partial sums with indices $n_k\sim cq^k$,
$1\le c<q$, converge weakly to the $c$-th convolution power of the
limit distribution for $c=1$.

L\'evy's remark quoted above is repeated, without proof, by Doeblin
\cite{doe1940} (see the paragraph after Th\'eor\`eme II on p.\ 78)
and with proof by Gnedenko (\cite{gn}, Theorem 3).
Gnedenko uses this fact to prove that if a distribution $G$ is not
stable, then the class $\mathcal{G}$ of different (modulo linear
transformations) convolution powers $G^{*c}$ is uncountable (see
also Gnedenko and Kolmogorov \cite{gk}, p.\ 189). Cs\"org\H{o}
(\cite{cs2}, Theorem 10) proved that $\cal{G}$ actually has the
cardinality of the continuum.

Note that L\'evy's proof of the limit relation $\hat{S}
\overset{d}{\longrightarrow} {\cal L}$ above, which expresses his
claim for $n=q^k$, is a coupling argument, defining the sequence
$(x_n)$ and the Poisson r.v.\ $N$ on the same probability space
and comparing pointwise the sequences  $\hat{S}=n^{-1/\alpha}
\sum_{k=1}^n x_k$ and $\bar{S}_t=n^{-1/\alpha} \sum_{k=1}^N
x_k$ ($t=n^{-1/\alpha}$), the latter of which converges weakly
to $\mathcal{L}$. Note the similarity of this method to Skorohod
embedding: L\'evy's method approximates St.\ Petersburg sums by
randomized i.i.d.\ sums, while Skorohod embedding represents
partial sums of square integrable i.i.d.\ sequences
as a randomly stopped Wiener process. In a sense,
L\'evy's idea is complementary to Skorohod's: it represents not
the partial sums of the $x_\nu$ themselves, but, after a small
perturbation, the limiting semistable variable. The method also
has important consequences for general i.i.d.\ sums. Let $x_n$ be
i.i.d.\ random variables concentrated on $[1, \infty)$ with tails
$G(x)=P(x_1> x)$, let $S_n=\sum_{k=1}^n x_k$, let $N$ be a Poisson
random variable with mean $n$, independent of the $x_k$'s, and let
$(a_n)$ be a positive numerical sequence. Clearly, $S_N/a_n$ has
the same distribution as the sum of points in a nonhomogeneous
Poisson process $\mathcal{P}$ in $(1/a_n, \infty)$ with intensity
measure with tails $T(x)=nG(a_nx)$. This distribution is
infinitely divisible with L\'evy measure concentrated on $(1/a_n,
\infty)$ with tails $nG(a_nx)$, and thus, apart from the fact that
the support of the L\'evy measure is $(1/a_n, \infty)$ instead of
$(0, \infty)$, it is the accompanying infinitely divisible
distribution to $S_n/a_n$ playing a central role in the Fourier
analytic theory (see e.g.\ \cite{gk}, p.\ 98). As before, under
suitable assumptions on $x_k$, the variables $S_n/a_n $ and
$S_N/a_n$ are close to each other pointwise and thus we get a
strong approximation of $S_n/a_n$ with an infinitely divisible
variable. In particular, it follows that the weak convergence of
$S_n/a_n$ along the whole sequence of integers or along a
subsequence is equivalent to the convergence of the corresponding
L\'evy measures $ndG(a_n x)$.
This is the well known 'generic' condition for weak convergence of i.i.d.\ sums,
used e.g.\ to characterize domains of attraction and domains of partial attraction.
It is typically proved by Fourier analytic methods and L\'evy's pointwise approximation
approach has considerable methodological interest.

For historical accuracy, one has to point out that strong approximation
is nowhere mentioned in L\'evy's paper and his sole interest was weak
convergence of i.i.d.\ sums. Neither is the St.\ Petersburg paradox
mentioned anywhere in \cite{levy35}, a curious fact, since L\'evy was
interested in the problem and discussed it in length in his book \cite{levy25},
pp.\ 122-133.\footnote{See Cs\"org\H{o} and Simons \cite{cssi}, pp.\ 69--70.}
The great power of strong approximation in probability theory
and statistics was not recognized until Strassen's paper \cite{str} and strong
approximation results based on the quantile transform such as those
of Koml\'os, Major and Tusn\'ady \cite{kmt},
or the strong approximation approach to domains of attraction and partial attraction due to
Cs\"org\H{o}, Haeusler and Mason \cite{cshm1} give much better rates
and have  a wider scope of applications than L\'evy's approach
above. We refer to \cite{csr},  \cite{cshm1}, \cite{sw} and
the references therein for history and further applications of the
quantile transformation method.

As far as the distributional closeness of i.i.d.\ sums and their
accompanying infinitely divisible laws is concerned, this has a
wide literature, starting with Doeblin \cite{doe1939}. L\'evy's
observation that the accompanying infinitely divisible
distribution of the partial sums $S_n=\sum_{k=1}^n x_k$ is the
same as the distribution of $S_N$ with a Poisson distributed $N$
has not gone unnoticed, see e.g.\ LeCam \cite{lecam}. However, the
usual path of estimation of the distributional closeness of $S_n$
and $S_N$ in the literature is not the direct path above, but,
following Doeblin \cite{doe1939}, it utilizes concentration
function arguments (also going back to L\'evy).
For the uniform distance of the distributions of $S_n$ and $S_N$
Kolmogorov \cite{ko56} obtained the bound $Cn^{-1/5}$ with an
absolute constant $C>0$; this has been improved by several
authors, see Arak and Zaitsev \cite{az} for a historical account
and optimal results.

\section{Asymptotics via the quantile transform}

In this section we discuss another argument in L\'evy \cite{levy35}, standing somewhat apart
from the main line of discussion of \cite{levy35}, but having profound consequences
on the behavior of i.i.d.\ sums. On pp.\ 372--374 of \cite{levy35} he makes the
following remark.

\bigskip\noindent
\begin{mypar}{0.7cm}{0.7cm}
{\bf Classification of the terms of  $S_n$ in decreasing order.}
{\it Let $\xi_1, \xi_2, \ldots, \xi_n$ denote the numbers $x_1,
\ldots, x_n$ ordered  in decreasing absolute values. The role
played by the largest of the $|x_\nu|$ in the previous discussions
makes one think that there is interest in studying the random
variables $\xi_p$ and considering $S_n$ as their sum.

$\phantom{81}$ Put $y_p=
F\{ |\xi_p|\}$. The numbers $y_1, y_2, \ldots, y_n$ are random
variables chosen at random between 0 and 1 and then arranged
in increasing order. Obviously,
$$ \mathcal{P} \{ y<y_p<y+dy\}=pC_n^p y^{p-1}(1-y)^{n-p}dy$$
from where we get easily, by using Euler integrals, that
$$ \mathcal{E} \{y_p\}=\frac{p}{n+1}, \quad \sigma^2 \{y_p\}=\frac{p(n-p+1)}{(n+1)^2 (n+2)}.$$
If $p\sim \alpha n$ as $n\to\infty$, the previous expressions
are equivalent to $\alpha$, resp.\
$\frac{\alpha(1-\alpha)}{n}$, as one can deduce it immediately
from Bernoulli's theorem. The case which is most interesting for us
is when $p$ is fixed and $n\to\infty$; then we have
$$\mathcal{E}
\{y_p\} \sim \frac{p}{n}, \quad \sigma^2\{y_p\} \sim \frac{p}{n^2} $$
\begin{equation}\label{gamma}
\mathcal{P} \{ ny_p < \eta\}\sim  \frac{1}{(p-1)!}\int_0^\eta \eta^{p-1} e^{-\eta} d\eta
\end{equation}
and thus, asymptotically, $ny_p$ is a random variable with mean
$p$ and mean quadratic deviation $\sqrt{p}$ and the normed difference $\frac{ny_p-p}{\sqrt{p}}$
is asymptotically normal (provided $p$ is very large). Further, between the
$y_p$ there is a positive correlation: if $y_p$ is known, then
$y_{p+q}$ can be considered as the $q$-th of $n-p$ variables,
chosen between $y_p$ and $1$ and arranged in increasing
order. The fact that each $y_p$ differs little from its
expectation permits us to neglect this correlation.

$\phantom{81}$
To fix the ideas, assume that the studied law is symmetric; each
$\xi_p$ is then  a variable with modulus $\rho_p$ known in terms
of the $y_p$ and with random sign. The sum $S_n$ can then be
written as
\begin{equation} \label{levydecomp}
\pm \varrho_1 \pm \varrho_2 + \ldots \pm \varrho_p \pm \ldots
\end{equation}
and here the number of terms grows with $n$, so that it can be
considered as a series.  If $\frac{\varrho_p}{\varrho_1}$ has an
order of magnitude independent of $n$,  $S_n$ has the order of
magnitude $\varrho_1$ or larger, according as the probability of the
convergence of the series is 1 or 0, i.e.\ according as the series $\sum
\varrho_k^2$ converges or not.  It is in the case of divergence of
this series when the law of large numbers applies.

$\phantom{8}$ Thus, if $F(x)\sim x^{-\alpha}$, we have (using $p/n$ as an approximative
value of $y_p$)
$$ \varrho_p^\alpha\sim \frac{n}{p}, \quad \frac{\varrho_p}{\varrho_1}
\sim p^{-\frac{1}{\alpha}}, $$
and thus the law of large numbers applies if $\sum
p^{-\frac{2}{\alpha}}$  diverges, that is if $\alpha\ge 2$ and
only in this case. We can thus reprove well known results by an
intuitive procedure which we found useful to mention, since the
validity of law of large numbers is connected with the divergence
of $\sum \varrho_p^2$. }
\end{mypar}

\bigskip
This remark seems to be the first application of the quantile
transform in probability theory, an argument used 40 years later
with spectacular success to the asymptotic study of sums of
independent random variables. To understand what is actually
stated and proved here, one has to note that the phrase "the law
of large numbers applies" for a sum of independent random variables
means, in L\'evy's terminology (see the top of p.\ 363 in
\cite{levy35}) something different, namely the uniform asymptotic negligibility
of the terms of the sum compared with the sum.
As L\'evy proves in the remark above, this condition holds for
sums of i.i.d.\ symmetric random variables with tails $\sim
cx^{-\alpha}$ iff $\alpha\ge 2$. But a more important consequence
of the argument above is  that for $0<\alpha<2$ the $p$-th largest
term of a sum of $n$ i.i.d.\ symmetric random variables with tails
$\sim cx^{-\alpha}$  has the order of magnitude
$(n/p)^{1/\alpha}$, i.e., the extremal terms of the sum have the
same order of magnitude as the sum itself. This is a fundamental
observation, the 'signature' property of i.i.d.\ sums with stable
tails. Naturally, $\rho_p$ depends also on $n$ and relation
(\ref{gamma}) shows that for fixed $p$ and $n\to\infty$ the limit
distribution of $ny_p$  is Gamma\,$(p, 1)$ (a fact also discovered
much later).  Thus using $F^{-1}(x) \sim c_1 x^{-1/\alpha}$
$(x\to 0$) and the connection between $y_p$ and $|\xi_p|$ we see that
$n^{-1/\alpha} \rho_p$ converges weakly, for fixed $p$ and $n\to\infty$, to
$1/Z_p^{1/\alpha}$ where $Z_p$ is a  Gamma\,$(p, 1)$
random variable. It is now tempting to divide in (\ref{levydecomp})
by $n^{1/\alpha}$, let $n\to\infty$ and conclude that a random
variable having the stable limit distribution of $S_n/n^{1/\alpha}$  has the infinite
representation
\begin{equation}\label{lwz}
\pm Z_1^{-1/\alpha} \pm Z_2^{-1/\alpha} \pm \ldots
\end{equation}
where $Z_p$ is a Gamma\,$(p, 1)$ variable.
L\' evy's  justification of this transition to infinite series is a bit vague: {\it "here the
number of terms grows with $n$, so that it can be considered as a series"} and his proof of the
a.s.\ convergence of the infinite series is not complete either: for justifying the use of the Kolmogorov
two-series criterion  he remarks merely that the correlation between the terms of the sum can be neglected.
The decomposition (\ref{lwz})
was proved only in 1981 by LePage, Woodroofe and Zinn
\cite{lwz} with the crucial addition that the $Z_p$'s  are the
partial sums of a {\it single} i.i.d.\ sequence of exponential
random variables with mean 1; this yields naturally the precise
joint distribution of the terms of the expansion.
The modern tool used in \cite{lwz} is
the representation  of the uniform ordered sample in the form
$(Z_1/Z_{n+1}, \ldots, Z_n/Z_{n+1})$ (see e.g.\ \cite{sw}, p.\ 335) and the proof of the a.s.\
convergence of the sum in Le Page, Woodroofe and Zinn \cite{lwz}
requires a delicate argument. As we will show in \cite{buj}, this
decomposition, adapted to semistable variables, leads to a strong
approximation of St.\ Petersburg sums with a semistable L\'evy
process with a.s.\ error $O(n^{1/2+\ve})$ and an asymptotically
normal remainder term, proving an unexpected central limit theorem
in St.\ Petersburg theory.

\bigskip\noindent
{\bf Acknowledgement.} The author thanks Professor Mikl\'os
Cs\"org\H{o} for his valuable comments.

\end{document}